%% file: coeffzeros.arxiv.tex
\documentclass[12pt]{amsart}

\usepackage{amssymb,amsmath,latexsym} \usepackage{rotating}
\usepackage{graphicx}
\usepackage{color}
\usepackage[letterpaper,textwidth=6.5in,marginparwidth=.15in]{geometry}
\usepackage{paralist}
\usepackage{enumerate}

\newtheorem{theorem}{Theorem}[section] 
\newtheorem{lemma}[theorem]{Lemma}
\newtheorem{corollary}[theorem]{Corollary}
\newtheorem{proposition}[theorem]{Proposition}

\newtheorem{conjecture}[theorem]{Conjecture}

\theoremstyle{definition} 
\newtheorem*{remark}{Remark}

\DeclareMathOperator{\conv}{conv}

\newcommand\vol{\operatorname{vol}}
\newcommand\aff{\operatorname{aff}}

\def\R{\mathbb{R}} \def\C{\mathbb{C}} \def\Z{\mathbb{Z}}
\def\N{\mathbb{N}} 
 
\def\O{\mathcal{O}}
 
\newcommand\x{\mathbf{x}}
\newcommand\y{\mathbf{y}}  \def\x{{\bf x}}

\title{Coefficients and Roots of Ehrhart Polynomials}

\author[Beck, De Loera, Develin, Pfeifle, and Stanley]{M. Beck, J. A. De Loera, M. Develin, J. Pfeifle, and R. P. Stanley}  

\address{Max-Planck-Institut f\"ur Mathematik, Vivatsgasse 7, 53111 Bonn, Germany}
\email{beck@mpim-bonn.mpg.de}

\address{Department of Mathematics, University of California, One Shields Avenue, Davis, CA 95616-8633, USA}
\email{deloera@math.ucdavis.edu}

\address{Department of Mathematics, University of California, Berkeley, California 94720, USA}
\email{develin@post.harvard.edu}

\address{Institut de Matem\`atica, Universitat de Barcelona, Gran Via de les Corts Catalanes 585, E-08007 Barcelona, Spain}
\email{julian@imub.ub.es}

\address{Department of Mathematics 2-375, Massachusetts Institute of Technology, Cambridge, MA 02139-4307, USA}
\email{rstan@math.mit.edu}

\begin{document}

\def\currentvolume{} \def\currentissue{} \pagespan{1}{60} \PII{}
\copyrightinfo{}{} \keywords{} \subjclass[2000]{}
\setlength{\parindent}{0pt}

\abstract The Ehrhart polynomial of a convex lattice polytope counts
integer points in integral dilates of the polytope. We present new
linear inequalities satisfied by the coefficients of Ehrhart
polynomials and relate them to known inequalities. We also
investigate the roots of Ehrhart polynomials. We prove that for
fixed $d$, there exists a bounded region of $\C$ containing all roots of
Ehrhart polynomials of $d$-polytopes, and that all real roots of these
polynomials lie in $[-d, \lfloor d/2 \rfloor)$. In contrast, we prove that when the
dimension $d$ is not fixed the positive real roots can be arbitrarily large.  
We finish with an experimental investigation of the Ehrhart polynomials
of cyclic polytopes and $0/1$-polytopes. \endabstract

\maketitle

\setlength{\parskip}{0.4cm} \bibliographystyle{amsplain}


\section{Introduction} 

In this article, a \emph{lattice polytope} $P \subset \R^d$ is a
convex polytope whose vertices have integral coordinates.  (For all
notions regarding convex polytopes we refer to~\cite{ziegler}.)  In
1967 Eug\`ene Ehrhart  proved that the function which
counts the lattice points in the $n$-fold dilated copy of $P$,
\[ 
    i_P: \N\to\N, \qquad i_P(n)\ = \ \# \left( n P \cap \Z^d \right)\,, 
\] 
is a polynomial in~$n$ (see \cite{ehrhart1,ehrhart2} and the
description in \cite{ehrhartbook}). In particular, $i_P$ can be
naturally extended to all complex numbers~$n$.  In this paper we
investigate linear inequalities satisfied by the coefficients of
Ehrhart polynomials and the distribution of the roots of Ehrhart
polynomials in the complex plane.

The coefficients of Ehrhart polynomials are very special. For example,
it is well known that the leading term of $i_P (n) $ equals the volume
of $P$, normalized with respect to the sublattice $\Z^d \cap \aff
(P)$. The second term of $ i_P (t) $ equals half the surface area of
$P$ normalized with respect to the sublattice on each facet of $P$,
and the constant term equals $1$. Moreover, the function $i_P^\circ
(n)$ counting the number of interior lattice points in~$nP$ satisfies
the \emph{reciprocity law} $ i_P (-n) =(-1)^{\dim P} i_P^\circ
(n)$~\cite{ehrhartbook,macdonald,stanleyreciprocity}.

Our first contribution is to establish new linear relations satisfied
by the coefficients of all Ehrhart polynomials. This is a continuation
of the pioneering work of Stanley, Betke \& McMullen, and Hibi
\cite{stanleydecomp,stanleyh1,BetkeMcMullen,hibi1}, who established
several families of linear inequalities for the coefficients (see
Theorems \ref{stanleylemma} and \ref{bmineqs}).  If we think of an
Ehrhart polynomial $i_P(n)=c_dx^d+c_{d-1}x^{d-1}+\dots+c_1x+1$ as a
point in $d$-space, given by the coefficient vector
$(c_d,c_{d-1},\dots,c_1)$, their results imply that the Ehrhart
polynomials of all $d$-polytopes lie in a certain polyhedral
complex. Betke and McMullen raised the issue \cite[page
262]{BetkeMcMullen} of whether other linear inequalities are possible. We
were indeed able to find such new inequalities
in the form of bounds for the \emph{$k$-th difference} of the Ehrhart
polynomial $i_P(n)$. These are defined recursively via
\[ 
    \ \Delta i_P (n) \ = \ i_P (n+1) - i_P (n)
\] 
and
\[ 
   \Delta^k i_P (n) \ = \ \Delta \left( \Delta^{ k-1 } i_P (n) \right)
   \quad\text{for }\ k\ge1 \quad\text{and} \quad \Delta^0 i_P(n)=i_P(n). \
\]
%
%
Our first result (proved in Section \ref{coeff}) is as follows.
\begin{theorem}\label{thm:diffinequ}
  
  If the lattice $d$-polytope $P\subset\R^d$ has Ehrhart polynomial
  $i_P (n) = c_d \, n^d + \dots + c_0$, then 
  \[ 
     {d\choose \ell}\Delta^k i_P (0)\ \leq \ \binom d k \Delta^\ell i_P (0)
     \qquad\text{for }\ 0\le k<\ell\le d.
  \] 
  In particular (put $k=0$ resp.\ $\ell=d$),
  \[ 
     \binom{d}{k} \ \le \ \Delta^k i_P (0)\ \le\ \binom d k d!  \, c_d 
     \qquad\text{for }\ 0\le k \le d. 
  \]

\end{theorem}

In Section \ref{coeff} we give a proof of Theorem
\ref{thm:diffinequ} using the language of rational generating
functions as established in~\cite{BetkeMcMullen,stanleyec1}, and
make a summary of known linear constraints and their strength.

The relation between the coefficients and the roots of
polynomials, via elementary symmetric functions, suggests that once
we understand the size of the coefficients of Ehrhart polynomials 
we could predict the distribution of their roots in the complex plane.
The second contribution of this paper is a general study
of the roots of Ehrhart polynomials. 

There is clearly something special about the roots of Ehrhart
polynomials.  Take for instance the integer roots: Since a lattice
polytope always contains some integer points (namely, its vertices),
all integer roots of its Ehrhart polynomial are negative. More
precisely, by the reciprocity law, the integer roots of an Ehrhart
polynomial are those~$-n$ for which the open polytope~$n P^\circ$
contains no lattice point.  For instance, the Ehrhart polynomial
$\binom{ n+d }{ d }$ of the \emph{standard simplex} in~$\R^d$ (with
vertices at the origin and the unit vectors on the coordinate axes)
has integer roots at $n=-d, -d+1, \dots, -1$.

The roots of the Ehrhart polynomial of the cross polytope
\[ 
  \O^d \ = \ \left\{ (x_1, \dots, x_d) \in \R^d : \ \left| x_1 \right| +
  \dots + \left| x_d \right| \leq 1 \right\} ,
\] 
also exhibit special behavior: Bump et al.
\cite{bumpchoikurlbergvaaler} and Rodriguez \cite{rodriguez} proved
that the zeros of~$i_{\O^d}$ all have real parts equal to~$-1/2$.

Using classical results from complex analysis and the linear
inequalities of Theorem~\ref{thm:inequalities}, we derive in Section
\ref{roots} the
following theorems:

\begin{theorem}\label{thm:rootbound}
\begin{compactenum}
  
\item[\upshape{(a)}]\label{thm1:allroots} The roots of Ehrhart
  polynomials of lattice $d$-polytopes are bounded in norm by
  $1+(d+1)!\,$.
  
\item[\upshape{(b)}] \label{thm1:realroots} All real roots of Ehrhart
  polynomials of $d$-dimensional lattice polytopes lie in the
  half-open interval $[-d, \lfloor d/2 \rfloor)$.

\end{compactenum}

\end{theorem}

The upper bound we present in Theorem \ref{thm:rootbound} (b) is not
tight. For example, in Proposition~\ref{prop:dim4}, we give a very
short self-contained proof of the fact that Ehrhart polynomials for
polytopes of dimension $d \leq 4$ have real roots in the interval
$[0,1)$. In contrast with the above theorem we can also prove the
following result.

\begin{theorem} \label{growth}
For any positive real number $t$ there exist an Ehrhart polynomial of
sufficiently large degree with a real root strictly larger than
$t$. In fact, for every $d$ there is a $d$-dimensional $0/1$-polytope 
whose Ehrhart polynomial has a real zero $\alpha_d$ such that
$\lim_{d\rightarrow \infty}\alpha_d/d = 1/(2\pi e) = 0.0585\cdots$.
\end{theorem}

Our third contribution is an experimental study of the roots and
coefficients of Ehrhart polynomials of concrete families of lattice
polytopes. Our investigations and conjectures are supported by
computer experimentation using {\tt LattE}~\cite{lattemanual,latte}
and {\tt polymake}~\cite{polymake}.  For the complex roots, we offer
the following conjecture, based on experimental data.

\begin{conjecture}
  All roots $\alpha$ of Ehrhart polynomials of lattice
  $d$-polytopes satisfy $-d \leq \text{\rm Re } \alpha \leq d-1$.
\end{conjecture}

We also computed the Ehrhart polynomials of all $0/1$-polytopes of
 dimension less than or equal to 4 and for many cyclic polytopes:

\begin{conjecture}\label{cyclicconj}
  For the cyclic polytope $C(n,d)$ realized with integral vertices on
  the moment curve $\nu_d (t) := \left( t, t^2, \dots, t^d \right) $,
  \[
     i_{C(n,d)}(m)=\vol(C(n,d)) \, m^d+i_{C(n,d-1)}(m).
  \]
  Equivalently,
  \[
     i_{C(n,d)}(m)= \sum_{ k=0 }^d \vol_k (C(n,k)) \, m^k.
  \]
\end{conjecture}

We have experimentally verified this conjecture in many cases.


\section{An appetizer: dimension two}

Since Ehrhart polynomials of lattice 1-polytopes (segments) are of the
form $\ell n+1$, where $\ell$ is the length of the segment, we know
everything about their coefficients and roots: the set of roots of
these polynomials is $\{-1/\ell:\ell\ge1\}\subset[-1,0)$.

The first interesting case is dimension $d=2$. Pick's Theorem tells 
us that the Ehrhart polynomial of a lattice 2-polytope $P$ is
\[
   i_P(n) = c_2 \, n^2 + c_1 \, n + 1 \ ,
\]
where $c_2$ is the area of $P$ and $c_1$ equals $1/2$ times the number 
of boundary integer points of $P$.
In 1976, Scott  established the following linear relations.
Two polytopes are \emph{unimodularly equivalent} if there is a function 
which maps one to the other and which preserves the integer lattice.

\begin{theorem} {\upshape \cite{scott}}
Let $ i_P(n) = c_2 \, n^2 + c_1 \, n + 1 $ be the Ehrhart polynomial 
of the lattice 2-polytope $P$. If $P$ contains an interior integer 
point, and $P$ is not unimodularly equivalent to 
$\conv \left\{ (0,0), (3,0), (0,3) \right\}$,
then
\[
   c_1 \ \le \ \frac 1 2 c_2 + 2 \ .
\]
\label{T:scott}
\end{theorem}
By Pick's Theorem, for 2-polytopes with
no interior lattice points, we have $c_1 = c_2 + 1$.
For $P = \conv \left\{ (0,0), (3,0), (0,3) \right\}$, we obtain
$i_P(n) = 9/2 \, n^2 + 9/2 \, n + 1$.

\begin{figure}[htbp]
  \centering 
  \includegraphics[height=9cm]{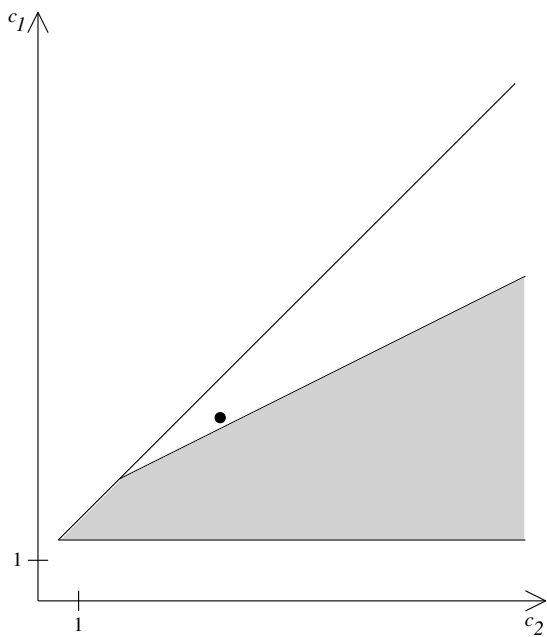}
  \caption{Regions in which Ehrhart polynomials of lattice 2-polytopes lie.
         It consists of 3 half lines, an open region (only points with half-integral coordinates are possible), plus an exceptional point.} 
  \label{fig:dim2cone}
\end{figure}

It is interesting to ask which degree-2 polynomials can possibly be
Ehrhart polynomials. Since the constant term has to be 1, we can think
of such a polynomial as a point $(c_2,c_1)$ in the plane. From the
geometry of lattice 2-polytopes, we know such an Ehrhart polynomial
must have half-integral coordinates. Aside from Scott's inequality, we
can trivially bound $c_1 \ge 3/2$, since every lattice 2-polytope has
at least 3 integral points, namely its vertices. From these
considerations, we arrive at Figure \ref{fig:dim2cone}, which shows
regions of possible Ehrhart polynomials of 2-polytopes.

Depicted are (parts of) three lines:
\begin{enumerate}[(i)]
\item\label{line1} $c_1 = 3/2$
\item\label{line2} $c_1 = c_2/2 + 2$
\item\label{line3} $c_1 = c_2 + 1$
\end{enumerate}
and the point $(c_2,c_1) = (9/2,9/2)$.  The ray (\ref{line1}) shows
the lower bound $c_1 \ge 3/2$.  This is a sharp lower bound, in the
sense that we can have polygons with exactly three boundary integer
points but arbitrarily large area.  The ray (\ref{line2}) is Scott's
bound, and the point $(c_2,c_1) = (9/2,9/2)$ corresponds to the
``exceptional'' polytope $\conv \left\{ (0,0), (3,0), (0,3) \right\}$
in Theorem~\ref{T:scott}.  The rectangles $\conv \left\{ (0,0), (2,0),
(2,x), (0,x) \right\}$, where $x$ is a positive integer, show that
there is a point on (\ref{line2}) for every half integer. Finally,
(\ref{line3}) corresponds to 2-polytopes which contain no interior
lattice point. There is a point on (\ref{line3}) for every half integer,
corresponding to the triangles $\conv \left\{ (0,0),
(1,0), (0,x) \right\}$ for a positive integer $x$.  The rays
(\ref{line1}) and (\ref{line3}) meet in the point $(1/2,3/2)$, which
corresponds to the standard triangle $\conv \left\{ (0,0), (1,0),
(0,1) \right\}$.  So the polyhedral complex containing all Ehrhart
vectors consists of the polyhedron bounded by (\ref{line1}),
(\ref{line2}), and (\ref{line3}) (shaded in Figure
\ref{fig:dim2cone}), plus the ray (\ref{line3}), plus the extra point
$(c_2,c_1) = (9/2,9/2)$. In fact, only points with half-integral
coordinates inside the complex are valid Ehrhart vectors. From these
constraints, we can locate possible roots of Ehrhart polynomials of
lattice 2-polytope fairly precisely.

\begin{theorem}\label{T:dim2roots}
The roots of the Ehrhart polynomial of any lattice 2-polytope are contained in
\[
   \left\{ -2, -1, - \frac 2 3 \right\} \cup \left\{ x+iy \in \C : \ - \frac 1 2 \le x < 0 , \, |y| \le \frac{ \sqrt{ 15 } }{ 6 } \right\} .
\]
\end{theorem}

\begin{proof}
We consider three cases, according to Scott's Theorem \ref{T:scott}.
First, if the lattice 2-polytope $P$ contains no interior lattice point then
$i_P(n) = A n^2 + (A+1) n + 1$ (by Pick's Theorem), where $A$ denotes
the area of $P$. The roots of $i_P$ are at $-1$ and $-1/A$.
Note that $A$ is half integral. 

The second case is the ``exceptional'' polytope 
$P = \conv \left\{ (0,0), (3,0), (0,3) \right\}$
whose Ehrhart polynomial $i_P(n) = 9/2 n^2 + 9/2 n + 1$ has roots $-2/3$ and $-1/2$.

This leaves, as the last case, 2-polytopes which contain an
interior lattice point and which are not unimodularly equivalent
to $\conv \left\{ (0,0), (3,0), (0,3) \right\}$. 
The corresponding Ehrhart polynomials $i_P(n)
= c_2 n^2 + c_1 n + 1$ satisfy the Scott inequality $c_1 \le c_2/2 +
2$.  Note that (because $P$ has an interior lattice point) the area of
$P$ satisfies $c_2 \ge 3/2$. We have two possibilities:

\noindent (A) The discriminant $c_1^2 - 4 c_2$ is negative. Then the
real part of a root of $i_P$ equals $- \frac{ c_1 }{ 2 c_2 }$ (which is negative). By
Pick's Theorem $c_1= c_2-I+1$ where $I$ is the number of interior
lattice points, that is, $-\frac{c_1}{2c_2} = -\frac12 - \frac{1-I}{2c_2}$. 
For fixed area $c_2$, this fraction is minimized
when $I$ is smallest possible, that is $I=1$.
The imaginary part of a root of $i_P$ is plus or minus
\[
\frac{ 1 }{ 2 c_2 } \sqrt{ 4 c_2 - c_1^2 } \ \le \ \frac{ 1 }{ 2 c_2 }
\sqrt{ 4 c_2 - \frac{ 9 }{ 4 } } \ = \ \sqrt{ \frac{ 1 }{ c_2 } - \left(
    \frac{ 3 }{ 4 c_2 } \right)^2 } \ ;
\]
here we used $c_1 \ge 3/2$. As a function in $c_2$, this upper bound
is decreasing for $c_2 \ge 1$.  Since $c_2 \ge 3/2$, we obtain as
an upper bound for the magnitude of the imaginary part of a root
\[
\sqrt{ \frac 2 3 - \left( \frac 1 2 \right)^2 } \ = \ \frac{ \sqrt{ 15
  } }{ 6 } \ .
\]

\noindent (B) The discriminant $c_1^2 - 4 c_2$ is nonnegative. 
Then the smaller root of $i_P$ is

\begin{eqnarray*}
- \frac{ c_1 }{ 2 c_2 } - \frac{ 1 }{ 2 c_2 } \sqrt{ c_1^2 - 4 c_2 }
&\ge& - \frac 1 4 - \frac{ 1 }{ c_2 } - \frac{ 1 }{ 2 c_2 } \sqrt{
  \left( \frac{ c_2 }{ 2 } + 2 \right)^2 - 4 c_2 } \\ &=& - \frac 1 4 -
\frac{ 1 }{ c_2 } - \frac{ 1 }{ 2 c_2 } \left( \frac{ c_2 }{ 2 } - 2
\right) \ = \ - \frac 1 2
\end{eqnarray*}
(Note that in this case $c_2 \ge 4$.) 

Finally, the larger root is negative, since all the coefficients of $i_P$ are positive.
\end{proof}


\section{Linear inequalities for the coefficients of Ehrhart polynomials}
\label{coeff}

In this section, we prove Theorem~\ref{thm:diffinequ}, which bounds
the ratio of the $k$-th and $\ell$-th differences of any Ehrhart
polynomial solely in terms of $d$, $k$, and $\ell$. It is perhaps
worth observing that most of our arguments are valid for a somewhat larger 
class of polynomials. To describe this class, we define the 
\emph{generating function} of the polynomial~$p$ as
\[
   S_p(x) = \sum_{ n \geq 0 } p(n) \, x^n \,.
\]
  It is well known (see, e.g.,
\cite[Chapter 4]{stanleyec1}) that, if $p$ is of degree $d$, then
$S_p$ is a rational function of the form
\begin{equation}\label{eq:ratfunc}
   S_p(x) = \frac{ f(x) }{ (1-x)^{ d+1 }  }  \,,
\end{equation}
where $f$ is a polynomial of degree at most $d$.  Most of our results
hold for polynomials~$p$ for which the numerator of $S_p$ has only
nonnegative coefficients. Ehrhart polynomials are a particular case, as
seen from the following theorem of Stanley.

\begin{theorem} {\upshape \cite[Theorem 2.1]{stanleydecomp}} 
  Suppose $P$ is a convex lattice polytope. Then the generating
  function $\sum_{ n \geq 0 } i_P(n) \, x^n $ can be written
  in the form of~\eqref{eq:ratfunc}, where $f(x)$ is a
  polynomial of degree at most $d$ with nonnegative integer
  coefficients.
\label{stanleylemma}
\end{theorem}

Another well-known (and easy-to-prove) fact about rational generating
functions (see, e.g., \cite[Chapter 4]{stanleyec1}) is the following.

\begin{lemma} \label{lem:aexpr}
  Suppose that $p\in\R[n]$ is a polynomial of degree $d$ with
  generating function $S_p(x) = (a_d x^d + \dots + a_1x +
  a_0)/(1-x)^{d+1}$. Then $p$ can be recovered as
  \begin{eqnarray}\label{diffidentity0}
    p(n) & = & \sum_{ j=0 }^{ d } a_j \binom{ d+n-j }{ d } \,. \\
    \noalign{More generally, we have the identity}\notag\\
    \label{diffidentity}
       \Delta^k p(n) & = & \sum_{ j=0 }^{ d } a_j \binom{ d+n-j }{ d-k } \qquad\text{for } k\ge0.
  \end{eqnarray}
\end{lemma}

\begin{proof} Equation \eqref{diffidentity0} follows from expanding $1/(1-x)^{ d+1 }$
  into a binomial series. For~\eqref{diffidentity}, we proceed by induction on
  $k$.  For $k=0$, the statement is (\ref{diffidentity0}), while for
  $k\ge1$ we have by the induction hypothesis
  \begin{align*}
    \Delta^k p(n)
    &\ = \ \Delta^{ k-1 } p(n+1) - \Delta^{ k-1 } p(n) \\
    &\ = \ \sum_{ j=0 }^{ d } a_j \left( \binom{ d+n+1-j }{ d-k+1 }-\binom{ d+n-j }{ d-k+1 } \right) \\
    &\ = \ \sum_{ j=0 }^{ d } a_j \binom{ d+n-j }{ d-k } \,.
  \end{align*}
\end{proof}

Combining Theorem \ref{stanleylemma} and Lemma \ref{lem:aexpr} immediately yields the following fact.

\begin{corollary}
  For any lattice polytope $P$ and $k\ge0$, we have
  $\Delta^k\,i_P(0)\ge0$.
\end{corollary}

\begin{proof} This follows because those binomial coefficients in the
final expression for $\Delta^k p(n)$ are either positive or zero.
\end{proof}

%

\begin{proof}[Proof of Theorem \ref{thm:diffinequ}]
  We will use the falling-power notation $d^{\underline{j}} =
  d(d-1)\cdots(d-j+1)$, along with the obvious relation
  $k^{\underline{j}}<\ell^{\underline{j}}$ for $j \leq k<\ell$, and the
  identity
  \[
    \binom{d-j}{d-k} \ = \ 
    \binom{d}{k}\frac{k^{\underline{j}}}{d^{\underline{j}}}.
  \]
  The statement now follows from Lemma \ref{lem:aexpr} (\ref{diffidentity}) by
  \[
    \binom{d}{\ell} \binom{d-j}{d-k} 
    \ = \ \binom{d}{\ell}\binom{d}{k}
    \frac{k^{\underline{j}}}{d^{\underline{j}}} 
    \ < \ \binom{d}{k}\binom{d}{\ell}\,
    \frac{\ell^{\underline{j}}}{d^{\underline{j}}} 
    \ = \ \binom{d}{k}\binom{d-j}{d-\ell}.
  \]
\end{proof}

Theorem~\ref{thm:diffinequ} is not the first set of linear
inequalities on coefficient vectors of Ehrhart polynomials. Indeed, in
1984, Betke and McMullen \cite[Theorem 6]{BetkeMcMullen} obtained the
following inequalities.

\begin{theorem} \label{bmineqs} 
  Let $P$ be a lattice $d$-polytope whose Ehrhart polynomial is
  $\sum_{i=0}^d c_i n^i$.  Then 
\[
   c_r \ \leq \ (-1)^{d-r} s(d,r)\, c_d + (-1)^{d-r-1}\frac{s(d,r+1)}{(d-1)!}
   \qquad\text{for }\ r = 1, 2, \dots, d-1,
\]
where $s(k,j)$ denote the Stirling numbers of the first kind. 
\hfill$\qed$
\end{theorem}

In that paper, Betke and McMullen sent out a challenge to the community to discover new inequalities 
for these coefficient vectors. The following theorem sums up the current state of affairs.

\begin{theorem} \label{thm:inequalities}
  Let $P$ be a $d$-dimensional lattice polytope, with Ehrhart
  polynomial $i_P(n) = \sum_{i=0}^d c_i n^i = \sum_{i=0}^d
  a_i\binom{n+d-i}{d}$. Then the following inequalities are valid for
  $0\le k<\ell\le d$ and $0\le i\le d$:
\begin{equation}\label{ineq:betke-mcmullen}
   c_r \ \leq \ (-1)^{d-r} s(d,r)\, c_d + (-1)^{d-r-1}\frac{s(d,r+1)}{(d-1)!} ,
\end{equation}
  \begin{center}
    \begin{minipage}{.48\linewidth}
      \begin{eqnarray}
        \label{ineq:diffrat}
        {d\choose k} \Delta^\ell\, i_P(0) &\ge& {d\choose \ell}
        \Delta^k\, i_P(0), \\
        \label{ineq:facet1}
        \binom{d+1}{2}\, c_d & \ge& c_{d-1},\\
        \label{ineq:facet2}
        i_P(1) &\ge& d+1,\\
        \label{ineq:iterated}
        \Delta^k\, i_P(0) &\ge& \binom d k ,
      \end{eqnarray}
    \end{minipage} \hfill
    \begin{minipage}{.48\linewidth}
      \begin{eqnarray} 
        \label{ineq:volume}
        c_d &\ge& c_0/d!,\\
        \label{ineq:surface}
        c_{d-1} &\ge& c_0\,\dfrac{d+1}{2(d-1)!},\\
        \label{ineq:interior}
        \sum_{i=0}^d (-1)^{d-i} c_i &\ge& 0,\\
        \label{ineq:apositive}
        a_i &\ge& 0 .
     \end{eqnarray}
   \end{minipage}
 \end{center}
Moreover,
\begin{align}
  a_d + a_{d-1} + \dots + a_{d-i} \quad & \le && \ a_0 + a_1 + \dots + a_i +
  a_{i+1}
  &&\text{for all }  \ 0\le i \le \lfloor (d-1)/2 \rfloor\,. \label{eq:hibi1}
  \\
\noalign{Whenever $a_s\ne0$ but $a_{s+1}=\dots=a_d=0$, then }
   a_0 + a_1 + \dots + a_i \quad & \le && \ a_s + a_{s-1} + \dots + a_{s-i}
   && \text{for all }\ 0\le i \le s\,; \label{eq:stanleyeq}\\
\noalign{finally, if $a_d\ne0$, then}
  a_1 \quad & \le && \ a_i &&\text{for all } \ 2\le i < d\,. \label{eq:hibi2}
\end{align}
\end{theorem}

\begin{proof}
  The inequalities \eqref{ineq:betke-mcmullen} and
  \eqref{ineq:diffrat} are the contents of Theorems~\ref{bmineqs} and
  \ref{thm:diffinequ}; while \eqref{ineq:facet1},
  \eqref{ineq:facet2}, and \eqref{ineq:iterated} are the special cases
  $(k,\ell)=(d-1,d)$, $(k,\ell)=(0,1)$, and $k=0$, respectively.
  \eqref{ineq:volume} and \eqref{ineq:surface} say that the volume and
  the normalized surface are at least as big as for a primitive
  simplex.  Inequality \eqref{ineq:interior} follows from Ehrhart
  reciprocity.  Inequality \eqref{ineq:apositive} is the statement of
  Theorem~\ref{stanleylemma}.  Incidentally, \eqref{ineq:volume} also
  follows from \eqref{ineq:iterated}, and \eqref{ineq:interior}
  follows from \eqref{ineq:apositive}, both by specializing to $i=d$.
  Inequality~\eqref{eq:stanleyeq} was proved by
  Stanley~\cite{stanleyh1}, and inequalities
  \eqref{eq:hibi1}~and~\eqref{eq:hibi2} by Hibi~\cite{hibi1,hibi2}.
\end{proof}

It is illuminating to compare these inequalities with each other.
Since inequality \eqref{ineq:apositive} was used to prove
Theorem~\ref{bmineqs} (by Betke and McMullen) and
Theorem~\ref{thm:diffinequ}, it seems stronger than the other
inequalities.  Indeed, the only inequality among
\eqref{ineq:betke-mcmullen}--\eqref{ineq:apositive} which does not
follow from \eqref{ineq:apositive} is \eqref{ineq:surface}.
Experimental data for small~$d$ shows that neither
\eqref{ineq:betke-mcmullen} nor \eqref{ineq:diffrat} imply
the other.

The set of linear inequalities of Theorem~\ref{thm:inequalities}
describes an unbounded complex of half-open polyhedra in $\R^{d+1}$
inside which all coefficient vectors of Ehrhart polynomials live.
From this, we obtain a bounded complex~${\mathcal Q}^d$ by cutting with
the normalizing hyperplane~$c_d = 1$. By~\eqref{eq:stanleyeq} each
constraint $a_s\ne0, a_{s+1}=\dots=a_d=0$ for $s=1,2,\dots,d$ defines
a half-open polytope~$E_s\in{\mathcal Q}^d$ of dimension~$s$ that is
missing one facet; $E_0$~is a single point.

\begin{figure}[htbp]
  \centering
  \includegraphics[width=.6\linewidth]{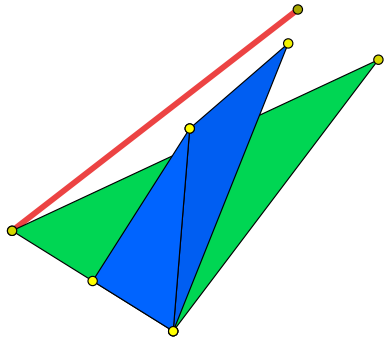}
  \caption{The complex ${\mathcal Q}^3$ of half-open polytopes, inside
    which the possible Ehrhart coefficients of all $3$-dimensional
    polytopes lie. The facets of the tetrahedron
    corresponding to $a_3\ne0$ are $a_0\ge0$, $c_2\ge1$,
    \eqref{eq:stanleyeq}~and~\eqref{eq:hibi2}; those of the triangle
    corresponding to $a_3=0$, $a_2\ne0$ are $a_0,a_1\ge 0$
    and~\eqref{eq:stanleyeq}; and those of the segment $a_2=a_3=0$,
    $a_1\ne0$ are $a_0\ge0$ and~\eqref{eq:stanleyeq}.}
  \label{fig:Q3}
\end{figure}

\begin{figure}[htbp]
  \centering
  \includegraphics[width=.6\linewidth]{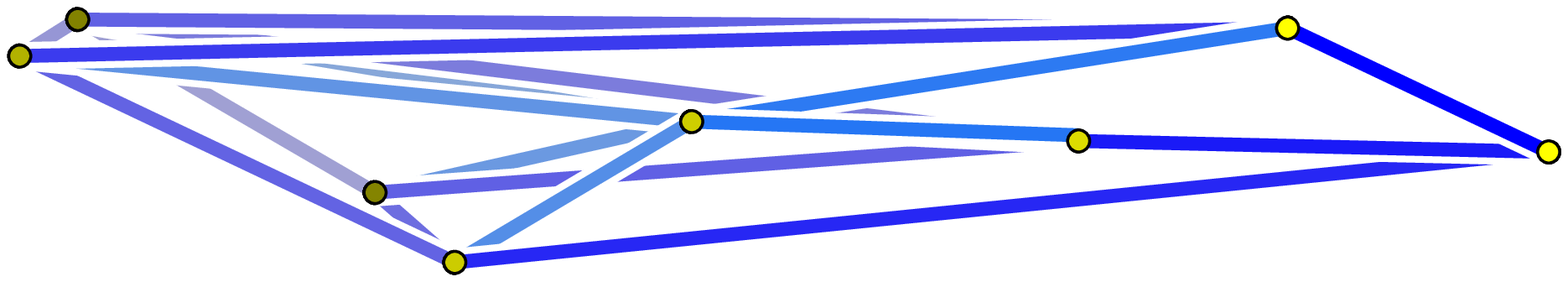}\\[-12ex]
  \includegraphics[width=.6\linewidth]{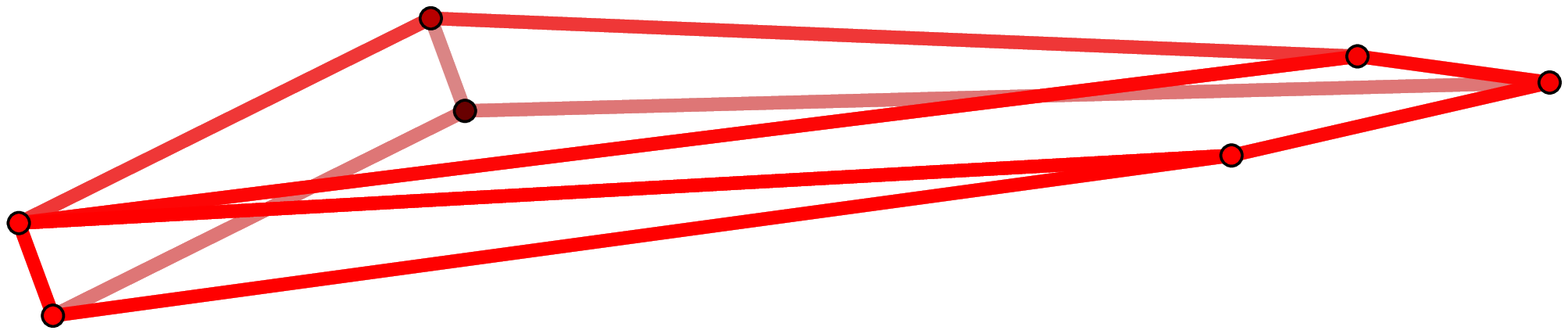}
  \caption{The $4$-dimensional (top) and $3$-dimensional (bottom)
    member of the complex ${\mathcal Q}^4$. The facets of the
    $4$-dimensional polytope are $a_0\ge0$, $c_3\ge 5/12$,
    \eqref{eq:stanleyeq}~for~$i=0$, \eqref{eq:hibi1}~for~$i=0,1$, and
    \eqref{eq:hibi2}~for~$i=2,3$; those of the $3$-dimensional one
    corresponding to $a_4=0$~but~$a_3\ne0$ are $a_0,a_1\ge0$, $c_3\ge
    5/12$, \eqref{eq:hibi1}~for~$i=1$, \eqref{eq:stanleyeq}~for~$i=1$,
    and \eqref{eq:hibi2}~for~$i=2$; etc.}
  \label{fig:Q4}
\end{figure}

Here are some particular cases: The bounded complex ${\mathcal Q}^3$
consists of one half-open $s$-dimensional simplex for each $s=0,1,2,3$
(Figure~\ref{fig:Q3}), and the half-open $3$- and $4$-dimensional
polytopes of~${\mathcal Q}^4$ are shown in Figure~\ref{fig:Q4}.

An important question about any linear inequality is whether or not it
defines a facet of ${\mathcal Q}^d$. We rephrase Betke and McMullen's
question \cite{BetkeMcMullen}:

{\bf Problem.}  \emph{ Are there other linear inequalities for the
coefficients of an Ehrhart polynomial aside from those in
Theorem~\ref{thm:inequalities}? Do they define facets of the
polyhedral complex inside which all coefficient vectors of Ehrhart
polynomials live?}


\section{The roots of Ehrhart polynomials}\label{roots}

When one has a family of polynomials, a natural thing to look at are
its roots. What is the general behavior of complex roots of Ehrhart
polynomials? As a consequence of the inequalities on its coefficients, 
we give bounds on the norm of roots of any Ehrhart polynomial in dimension $d$. 
The basis $\{ \binom{d+n-j}{d} : 0\le j\le d\}$ of the vector space of
  polynomials of degree~$d$ turns out to be much more natural than the
  basis $\{n^i:0\le i \le d\}$ for deriving bounds on the roots of
  Ehrhart polynomials $i_P(n)=\sum_{i=0}^d a_i\binom{n+d-i}{d} =
  \sum_{i=0}^d c_i n^i$. Also, recall the following classical result of Cauchy (see, for example,
\mbox{\cite[Chapter VII]{marden}}).

\begin{lemma}\label{mardenlemma}
  The roots of the polynomial $p(n) = c_d n^d + c_{ d-1 } n^{ d-1 } +
  \dots + c_0$ lie in the open disc
  \[ 
     \left\{ z \in \C : \ |z| < 1 + \max_{ 0 \leq j \leq d } 
       \left| \frac{ c_j }{ c_d } \right| \right\} . 
  \]
  \hfill $\qed$
\end{lemma}

Now we study roots of Ehrhart polynomials in general dimension. We
first give an easy proof bounding the norm of all roots.

\begin{proof}[Proof of Theorem~\ref{thm:rootbound}(a)]
  By Lemma \ref{mardenlemma} and Theorem \ref{bmineqs}, the maximal norm
  of the roots of $i_P$ is bounded by
\begin{eqnarray*}
  1 + \max_{ 0 \leq j \leq d } \left| \frac{ c_j }{ c_d } \right|
  &\le& 
  1+\max_{0 \leq j \leq d} \left| (-1)^{d-j} s(d,j) + 
    (-1)^{d-j-1} \frac{s(d,j+1)}{(d-1)!\, c_d} \right|  \\
  &\leq&  1+d!+d! d \ =\  1+(d+1)! \,.
\end{eqnarray*}
Here we have used the estimate $s(d,j) \leq |s(d,j)| \leq d!$
and the fact that $c_d \geq 1/d!$.
\end{proof}

While using crude estimates gives us a bound of $1+(d+1)!$, which
makes the main point that there exists a bound dependent only on $d$,
the actual bound on the roots can be improved greatly for specific
values of $d$. First of all, for small $d$, we can compute the inequalities exactly; here the 
inequalities from Theorem~\ref{thm:diffinequ} are used along with the Betke-McMullen inequalities. 
This gives appropriate bounds on the ratios of the coefficients of the Ehrhart polynomial. Second of 
all, Lemma~\ref{mardenlemma} is not the
best tool to use for specific cases, since calculating the
inequalities for small $d$ yields much lower bounds for $c_i/c_d$ when
$i$ is large. Instead, we use the following proposition.

\begin{proposition}[Theorem 27.1 \cite{marden}]
Let $p(n) = c_dn^d+c_{d-1}n^{d-1}+\cdots+c_0$ be a polynomial. Then the maximal value of the norm 
of a root of $p(n)$ is the value of the maximal root of $p^\prime(n) = |c_d| 
n^d-|c_{d-1}|n^{d-1}-|c_{d-2}|n^{d-2} -\cdots -|c_0|$.\hfill$\qed$
\end{proposition}

We use this and the exact calculation of the inequalities in question to obtain the following tighter
bounds on the roots of Ehrhart polynomials of $d$-polytopes.

\[
\begin{array}{r||r|r|r|r|r|r|r|r}
d & 2 & 3 & 4 & 5 & 6 & 7 & 8 & 9 \\ 
\hline
\text{bound} & 3.6 & 8.5 & 15.8 & 25.7 & 38.3 & 53.5 & 71.4 & 92.0
\end{array}
\]

The bound appears to grow roughly quadratically. We suspect
that there is a bound for the roots of Ehrhart polynomials of
$d$-polytopes which is polynomial in $d$. For real roots this is certainly the case;  we prove next that all real
roots of Ehrhart polynomials of $d$-polytopes lie in the interval
$[-d, \lfloor d/2 \rfloor)$. For this, we will use the following
well-known bound.

\begin{lemma} \emph{(Newton Bound)} Let $f\in\R[n]$ be a polynomial 
  of degree~$d$ and $B\in\R$ be such that all derivatives of~$f$ are
  positive at~$B$: $f^{(\ell)}(B)>0$ for $\ell=0,1,\dots,d$. Then all
  real roots of~$f$ are contained in $(-\infty, B)$. \hfill$\qed$
\end{lemma}

\begin{proof}[Proof of Theorem~\ref{thm:rootbound}(b)]
The lower bound follows from Theorem \ref{stanleylemma} and 
the simple observation that for (real numbers) $n < -d$ 
the binomial coefficients in 
\[
   i_P (n) = \sum_{ i=0 }^{ d } a_i \binom{ n+d-i }{ d }
\]
are all positive or all negative, depending on the parity of $d$. 

As for the upper bound, let $B=\lfloor d/2 \rfloor$. We now show that
$\alpha<B$ for any real root~$\alpha$ of~$i_P(n)$. For this, we will
make use of the fact that the second highest coefficient of any
Ehrhart polynomial measures half the normalized surface area. This
coefficient reads
\[
   c_{d-1} \ = \ \frac{1}{(d-1)!}\,\sum_{i=0}^d\, a_i\;(d-2i+1) 
\]
when expressed in terms of the $a_i$'s, so that the following
inequality is valid:
\begin{equation}\label{eq:surf}
   (d-1)!\, c_{d-1} \ = \  \sum_{i=0}^d\, a_i\;(d-2i+1) \ > \ 0 \,.
\end{equation}
Note that the coefficient~$s(i)=d-2i+1$ of $a_i$ in~\eqref{eq:surf} is
positive for $0\le i \le B$ and non-positive for $B+1\le i \le d$. We
now express the $\ell$-th derivative of $i_P$ evaluated at $n=B$ as
$i_P^{(\ell)}(B) = (\ell!/d!)\,\sum_{i=0}^d a_i\, g_i(B,\ell)$, and
claim that for  $0\le \ell\le d$, there exists a $\lambda(\ell)>0$
with
\[
    g_i(B,\ell) > \lambda(\ell)\, s(i)
    \qquad\text{for all }\ 0\le i \le d.
\]
This claim is the statement of Lemma~\ref{lem:lambda-ineq} below.  The
proof of Theorem \ref{thm:rootbound}(b) now follows from this
relation, inequality \eqref{eq:surf}, $a_0=1$ and $a_i\ge0$ for $1\le
i\le d$ via the following chain of inequalities:
\begin{eqnarray*}
  0 & < & \sum_{i=0}^d \big( g_i(B,\ell) - \lambda(\ell) \, s(i)\big) a_i \\
    & < & \sum_{i=0}^d \big( g_i(B,\ell) - \lambda(\ell) \, s(i)\big) a_i 
          + \lambda(\ell)\,\sum_{i=0}^d s(i)\,a_i\\
    & = & \sum_{i=0}^d g_i(B,\ell)\, a_i \  =\ i_P^{(\ell)}(B)\,.
\end{eqnarray*}
\end{proof}

\begin{remark}
It is a well-known fact that Ehrhart polynomials of lattice polytopes
form a special class of \emph{Hilbert polynomials}. More strongly,
they are special examples of Hilbert polynomials of Cohen-Macaulay
semi-standard graded $k$-algebras~\cite{stanleyh1} (this is
essentially the content of Theorem~\ref{stanleylemma}). It is then
natural to ask whether Ehrhart polynomials are special or whether the
bounds proved above hold in more generality. We stress that
inequality~\eqref{eq:surf}, used in previous arguments, comes from
geometric information about Ehrhart polynomials~$i_P(n)$. Indeed, from
the following proposition and Theorem \ref{thm:rootbound}(b), Ehrhart
polynomials are special in their root distribution:
\end{remark}

\begin{proposition} 
 For degree $d$ Hilbert polynomials associated to arbitrary
  semi-standard graded $k$-algebras the negative real roots are
  arbitrarily small and $d-1$ may appear as a root. In contrast, for
  fixed degree $d$, Hilbert polynomials of Cohen-Macaulay
  semi-standard graded $k$-algebras have all its real roots in the
  interval $[-d,d-1)$.
\end{proposition}

\begin{proof}
Indeed, it follows from \cite[Theorem 3.8]{Brenti} that for fixed $d$
and positive integers $a_0,\dots,a_d$, the polynomial
$a_0(x+a_1)(x+a_2)\dots(x+a_d)$ is the Hilbert polynomial of a
semi-standard graded $k$-algebra. Also, observe that the chromatic
polynomial of the complete graph on $d$ vertices has highest root
$d-1$, and that chromatic polynomials are known to be Hilbert
polynomials of standard graded algebras by a result attributed to
Almkvist (see the proof given by Steingr\'{\i}msson
~\cite{steingrimsson}).  Thus the first statement holds.

Now, in a Cohen-Macaulay semi-standard graded algebra, the Hilbert
polynomial can be written as $p(n)=\sum_{i=0}^d a_i\binom{n+d-i}{d}$,
where $a_i\ge0$ for $0\le i\le d$. Observe that all the binomial
coefficients in $p(n)$ are positive for (real numbers) $n > d-1$,
which establishes the upper bound of $d-1$. For the lower bound,
observe that for (real numbers) $n < -d$ all the binomial coefficients
are positive, respectively negative, depending on the parity of
$d$.
\end{proof}

To complete the proof of Theorem 1.2 (b), we need only to prove the
following lemma.

\begin{lemma} \label{lem:lambda-ineq}
  Fix $0\le\ell\le d-1$ and consider again the functions
  $s,g:\{0,1,\dots,d\} \to \Z$ defined by $s(i)=d-2i+1$ and
  $g(i)=g_i(B,\ell)$. Moreover, if we set
\begin{equation}\label{eq:lambda}
   \lambda(\ell)\ = \ \frac{1}{2}\, \big(g(B) - g(B+1)\big) \ = \ 
      \frac{d}{2}\,\sum_{I\in\binom{[d-1]}{d-\ell-1}}
      \prod_{k\in I} \, (d-k) \ > \ 0\,,
\end{equation}
then 
\begin{equation}
  \label{eq:final}
  g(i) \ \ge\ \lambda(\ell)\, s(i)
  \qquad\text{for }\ i=0,1,\dots, d\, .
\end{equation}
\end{lemma}

For this, we will also
need to prove Lemma~\ref{lem:gi-ineqs} below.  We will write
$[d-1]_0=\{0,1,\dots,d-1\}$, $[d-1]=\{1,2,\dots,d-1\}$, and
$\binom{S}{t}$ for the set of all $t$-element subsets of the finite
set~$S$. Now we express $i_P(n)$~as
\[
   i_P(n) \ = \ \frac{1}{d!}\, \sum_{i=0}^d a_i
                \prod_{k=0}^{d-1}\, (n+d-i-k),
\]
so that the $\ell$-th derivative of~$i_P$ is
\begin{eqnarray*}
   i_P^{(\ell)}(n) & = & \frac{\ell!}{d!}\, \sum_{i=0}^d a_i
                         \sum_{I\in\binom{[d-1]_0}{\ell}}
                         \prod_{k\in[d-1]_0\setminus I} (n+d-i-k) \\
                   & = & \frac{\ell!}{d!}\, \sum_{i=0}^d a_i
                         \sum_{I\in\binom{[d-1]_0}{d-\ell}}     
                         \prod_{k\in I}\, (n+d-i-k) \,.
\end{eqnarray*}

Note that we now have an explicit formula for the coefficient of $a_i$
in $(d!/\ell!)\,i_P^{(\ell)}$:
\begin{equation} \label{eq:gi}
    g_i(n,\ell) \ = \ \sum_{I\in\binom{[d-1]_0}{d-\ell}}     
                         \prod_{k\in I}\, (n+d-i-k) \,. 
\end{equation}
The following lemma shows that the piece-wise linear function
interpolating $g:\{0,1,\dots,d\}\to\Z$, $g(i) = g_i(B,\ell)$ is
positive, and its slope weakly increases in the range $0\le\ell\le d$
and $0\le i\le B+1$.  See Figure~\ref{fig:g-graphs}.

\begin{lemma} \label{lem:gi-ineqs}
The following inequalities are satisfied for $0\le\ell\le d$ and $0\le i\le B+1$:
\begin{eqnarray}
  \label{eq:eq1}
  g_i(B, \ell) & > & 0, \\
  \label{eq:eq2}
  g_i(B, \ell) - g_{i+1}(B, \ell) & > & g_{i+1}(B, \ell) - g_{i+2}(B, \ell). 
\end{eqnarray}
 
\end{lemma}

\begin{proof}
  Equation~\eqref{eq:eq1} follows because $k\le d-1$ and $i\le B+1$
  imply $B+d-i-k\ge0$. To show~\eqref{eq:eq2}, we abbreviate
  $m:=B+d-i$ and inspect the difference
  \[
     g_i(B, \ell) - g_{i+1}(B, \ell) \ = \ 
     \sum_{I\in\binom{[d-1]_0}{d-\ell}} \prod_{k\in I}\, (m-k) \; -
     \sum_{J\in\binom{[d-1]_0}{d-\ell}} \prod_{k\in J}\, (m-k-1)\,.
  \]
  If $0\notin I$, then the term corresponding to~$I$ in the first sum
  cancels with the term corresponding to $J=\{i-1:i\in I\}$ in the
  second sum:
  \[
      \prod_{k\in I}\,(m-k) - \prod_{k\in J}\,(m-k-1) \ = \ \prod_{k\in
        I}\,\left((m-k) - \big(m-(k-1)-1\big)\right) \ = \ 0,
  \]
  so we are left with summing over the
  sets~$I\in\binom{[d-1]_0}{d-\ell}$ that contain~$0$ and the sets~$J$
  that contain~$d-1$.  But for such summation sets, the difference
  simplifies to
  \begin{eqnarray*}
     g_i(B, \ell) - g_{i+1}(B, \ell) & = & 
        (m-0)\!\!\sum_{I\in\binom{[d-1]}{d-\ell-1}} \prod_{k\in I}\,(m-k)
        \;-\; (m-d)\!\!\sum_{J\in\binom{[d-2]_0}{d-\ell-1}} 
        \prod_{k\in I}\,\left(m-(k+1)\right) \\
     & = & d\!\! \sum_{I\in\binom{[d-1]}{d-\ell-1}}
        \prod_{k\in I}\, (B-d-i-k)\,,
  \end{eqnarray*}
  and \eqref{eq:eq2} follows by comparing the expressions $g_i(B, \ell)
  - g_{i+1}(B, \ell)$ and $g_{i+1}(B, \ell) - g_{i+2}(B, \ell)$ term by
  term.
\end{proof}

In the following, we will use Iverson's notation (see
\cite{concretemath}): the expression $[S]$ evaluates to $1$ resp.\ $0$
according to the truth or falsity of the logical statement~$S$.

\begin{proof}[Proof of Lemma \ref{lem:lambda-ineq}]
  First note that $s(B)=1$ for even~$d$, so that 
  \[
     g(B) - \lambda s(B)\  = \ \tfrac12 g(B) + \tfrac12 g(B+1) \ > \ 0\,;
  \]
  for odd~$d$, we have $s(B+1)=0$.  Now note that the graph of (the
  piecewise-linear function interpolating)~$\lambda s$ is a line,
  while $g(B+1)>0$ by~\eqref{eq:eq1} and the slope of the graph of~$g$
  is weakly increasing on $[0,B+2]$ by~\eqref{eq:eq2} (see
  Figure~\ref{fig:g-graphs}); this proves~\eqref{eq:final} for $0\le
  i\le B+[d\text{ odd}]$.

\begin{figure}[htbp]
  \centering
  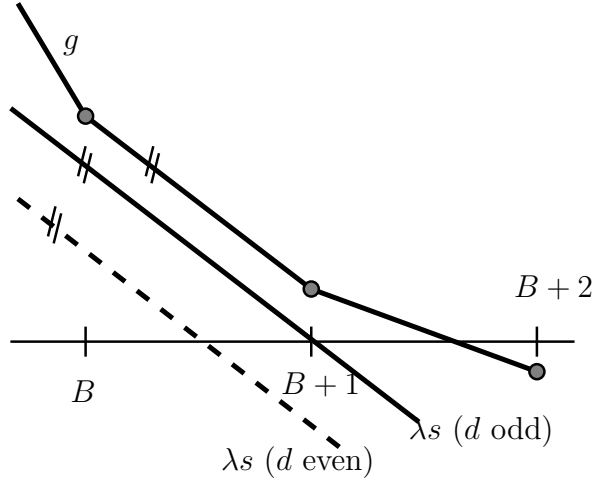
  \caption{The graphs of the functions $g$ and $\lambda s$ (solid for odd~$d$, 
    dashed for even~$d$). }
  \label{fig:g-graphs}
\end{figure}

Set $j=i-B$, so that we still need to prove~\eqref{eq:final} for
$1+[d\text{ odd}]\le j\le d-B$. By plugging
\eqref{eq:gi}~and~\eqref{eq:lambda} into~\eqref{eq:final} and
rearranging, we must show that for these values of $j$
\begin{equation} \label{eq:tailcase}
   \sum_{I\in\binom{[d-1]_0}{d-\ell}} \prod_{k\in I} (d-j-k) \;+\;
   \frac{d}{2}\,\big(2j-[d\text{ odd}]-1\big) \!\!
   \sum_{J\in\binom{[d-1]}{d-\ell-1}} \prod_{k\in J}(d-k) \ > \ 0\,.
\end{equation}
Note that each term in the second sum of \eqref{eq:tailcase} is positive, and
decompose the index sets~$I$ in the first sum into disjoint unions
$I=I_+\cup K$ such that $I_+\subset\{0,1,\dots,d-2j\}$ and
$K\subset\{d-2j+1,\dots,d-1\}$, and therefore $d-j-k>0$ for all $k\in
I_+$.
\begin{center}
  \begin{tabular}{c||cccccccccc}
    $k$ &              $0$ & $1$ & $\cdots$ & $d-2j$ &\vrule& $d-2j+1$ & $\cdots$ &
        $d-j$ & $\cdots$ & $d-1$
    \\\hline
    value of $d-j-k$ & $d-j$ & & $\cdots$ & $j$ &\vrule& $j-1$ & $\cdots$ &
        $0$ & $\cdots$ & $-(j-1)$\\\hline
    set in $I=I_+\cup K$
        & \multicolumn{4}{c}{$I_+$} &\vrule & \multicolumn{5}{c}{$K$}    
  \end{tabular}
\end{center}
If $|K|$~is odd, then the summand $\sigma(K)$ corresponding
to~$I_+\cup K$ cancels with the one corresponding to~$I_+\cup(d-j-K)$,
so we only need to consider even~$|K|$. In that case, $\sigma(K)>0$
(resp.\ $\sigma(K)<0$) if $|K\cap[d-j+1, d-1]|$ is even (resp.\ odd).  In
total, there are more than enough positive terms in~\eqref{eq:tailcase}
to cancel the negative summands.
\end{proof}

\begin{proposition}\label{prop:dim4}
  We have $\alpha<1$ for any real root $\alpha$ of an Ehrhart
  polynomial~$i_P$ of a lattice polytope~$P$ of dimension $d \leq 4$ .
\end{proposition}

\begin{proof}
It is enough to prove the statement in dimension 4 because of Theorem
\ref{thm:rootbound}(b). Suppose $f(n)=
p{n}^{4}+q{n}^{3}+r{n}^{2}+sn+1$ is the Ehrhart polynomial of a
lattice 4-polytope~$P$.  We know $p>0$ and $q>0$.  Because $f(1)$
counts the lattice points in $P$, we know that $p + q + r + s + 1 \geq
5$. By the reciprocity law, $f(-1)\geq 0$,
so $p - q + r - s + 1 \geq 0$. The top two coefficients of the shifted
polynomial $g(n)=f(n+1)=p n^4 + (4p+q)n^3 + g_2n^2 + g_1n + g_0$ are
positive, as is the constant term $g_0 = g(0) = f(1)$.  We will show
that $g_2$ and $g_1$ are nonnegative, and hence, by Descartes' rule of
signs, $g$ does not have a positive root. This implies that $f(n) =
g(n-1)$ does not have a real root larger than 1.  To prove that $g_2
\geq 0$, we add the inequalities $f(1) \geq 5$ and $f(-1) \geq 0$ to
obtain $2p + 2r \geq 3$ or $r \geq \frac 3 2 - p $, whence $g_2 = 6p +
3q + r \geq 5p + 3q + \frac 3 2 \geq 0$ (because $p,q \geq 0$).  A
similar reasoning yields \[ g_1 \ = \ 4p + 3q + 2r + s \ = \ (p+q+r+s)
+ (3p+2q+r) \ \geq\ 4 + 2p + 2q + \frac 3 2 \ \geq \ 0 \ ; \] here we
used the inequality $f(1) \geq 5$ again.  \end{proof}

We now conclude with the proof of Theorem \ref{growth}:

\begin{proof}[Proof of Theorem \ref{growth}] 
Given a positive integer $d$, consider the convex polytope $P_d$ defined
by the facet inequalities:
\[
0 \leq x_0 \leq x_k \leq 1 \quad \hbox{for} \quad 1 \leq k \leq d-1.
\]
$P_d$ is an \emph{order polytope} in the sense of
\cite{stanleyposetpolytopes} and thus it has $0/1$ vertices. We claim
that the Ehrhart polynomial of $P_d$ is given by $ i_{P_d}(n) =
(B_d(n+2) - B_d(0))/d$ where $B_d(x)$ is the
$d$-th Bernoulli polynomial.

Indeed, from the facet-defining inequalities of $P_d$ one sees that
$i_{P_d}(n)$ is the number of $d$-tuples of nonnegative integers
$(a_0, a_1, \dots, a_{d-1})$ such that $a_0 \leq a_k \leq n$. If $a_0 = j$
then there are $n-j+1$ choices for each $a_k$ $(k>0)$. Hence
$i_{P_d}(n) = \sum_{j=0}^n (n-j+1)^{d-1} = \sum_{j=1}^{n+1} j^{d-1}.$ A
classical identity of Bernoulli that says $\sum_{k=0}^{n-1}
k^{d-1}=(B_d(n)-B_d(0))/d$. Thus we get $i_{P_d}(n) = (B_{d}(n+2) -
B_{d}(0))/d$, a polynomial of degree $d$. Note that when $d$ is odd
then $B_d(0)=0$.  Finally, the results of \cite{veselovward} imply
that the largest real zero of $B_d(n)$ is asymptotically $d/(2\pi
e)$. Therefore, as stated, as the degree $d$ grows, the Ehrhart
polynomial of $P_d$ has larger and larger real roots. It is
worth remarking that since $d$ is the degree of the Ehrhart polynomial
of $P_d$ it differs from the upper bound of $\lfloor d/2 \rfloor$ in Theorem
\ref{thm:rootbound} only by a constant factor.
\end{proof}


\section{Special families of polytopes}

\begin{figure}[htbp]
  \centering
    \includegraphics[height=8cm]{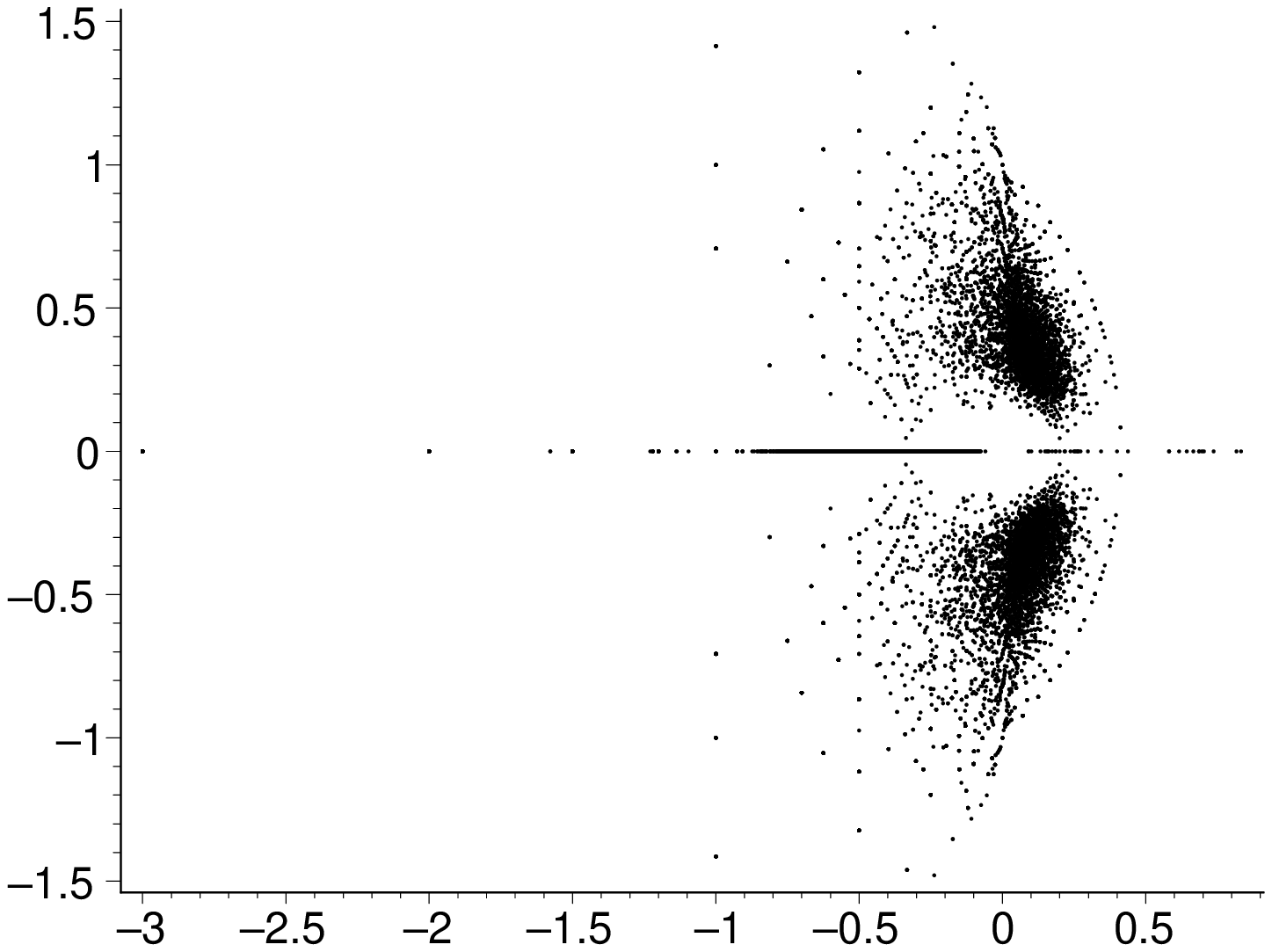} 
  \caption{
  The zeros of Ehrhart
  polynomials corresponding to $100000$~random $3$-dimensional lattice
  simplices.} \label{3dlatpoly}
\end{figure}

We begin this section with some charts showing the behavior of roots for
hundreds of Ehrhart polynomials computed using {\tt LattE} and {\tt
Polymake}. In Figure \ref{3dlatpoly} we show the distribution of roots
of a large sample of Ehrhart polynomials of lattice
3-polytopes. 

In Table \ref{Ziegler} we collected a small sample of Ehrhart polynomials of
$0/1$ polytopes and cyclic polytopes from the experiments we performed. 

{\tiny
\begin{sidewaystable}
  \centering
\begin{tabular}{|c|c|} 
\hline
name &  Ehrhart polynomial $P(s)$ \\ \hline 
cube & ${t}^{3}+3\,{t}^{2}+3\,t+1$ \\ \hline
cube minus corner & $5/6\,{t}^{3}+5/2\,{t}^{2}+8/3\,t+1$ \\ \hline
prism & $1/2\,{t}^{3}+2\,{t}^{2}+5/2\,t+1$ \\ \hline
nameless & $2/3\,{t}^{3}+2\,{t}^{2}+7/3\,t+1$ \\ \hline
octahedron & $2/3\,{t}^{3}+2\,{t}^{2}+7/3\,t+1$ \\ \hline
square pyramid & $1/3\,{t}^{3}+3/2\,{t}^{2}+{\frac {13}{6}}\,t+1$ \\ \hline
bypyramid  & $1/2\,{t}^{3}+3/2\,{t}^{2}+2\,t+1$ \\ \hline
unimodular tetrahedron & $1/6\,{t}^{3}+{t}^{2}+{\frac {11}{6}}\,t+1$ \\ \hline
fat tetrahedron  & $1/3\,{t}^{3}+{t}^{2}+5/3\,t+1$ \\ \hline
as:6-18.poly  &  ${\frac {83}{240}}\,{x}^{6}+{\frac {307}{240}}\,{x}^{5}+{\frac {41}{16}
}\,{x}^{4}+{\frac {217}{48}}\,{x}^{3}+{\frac {611}{120}}\,{x}^{2}+{
\frac {16}{5}}\,x+1 $  \\ \hline
cf:10-11.poly   &  ${\frac {11}{3628800}}\,{x}^{10}+{\frac {11}{725760}}\,{x}^{9}+{\frac {
17}{60480}}\,{x}^{8}+{\frac {121}{24192}}\,{x}^{7}+{\frac {7643}{
172800}}\,{x}^{6}+{\frac {8591}{34560}}\,{x}^{5}+{\frac {340873}{
362880}}\,{x}^{4}+{\frac {84095}{36288}}\,{x}^{3}+{\frac {59071}{16800
}}\,{x}^{2}+{\frac {7381}{2520}}\,x+1$ \\ \hline
cf:4-5.poly & $1/12\,{x}^{4}+1/2\,{x}^{3}+{\frac {17}{12}}\,{x}^{2}+2\,x+1$ \\ \hline
cf:9-10.poly & $ {\frac{1}{120960}}\,{x}^{9}+{\frac {1}{4480}}\,{x}^{8}+{\frac{61}{20160}}\,{x}^{7}+{\frac{79}{2880}}\,{x}^{6}+{\frac {997}{5760}}\,{x}^{5}+{\frac{4223}{5760}}\,{x}^{4}+{\frac{30043}{15120}}\,{x}^{3}+{\frac{32651}{10080}}\,{x}^{2}+{\frac{2383}{840}}\,x+1 $\\ \hline

cf:8-9.poly &${\frac {11}{40320}}\,{x}^{8}+{\frac {1}{1120}}\,{x}^{7}+{\frac {1}{64}}\,{x}^{6}+{\frac {9}{80}}\,{x}^{5}+{\frac {1039}{1920}}\,{x}^{4}+{\frac {267}
{160}}\,{x}^{3}+{\frac {5933}{2016}}\,{x}^{2}+{\frac {761}{280}}\,x+1$ \\ \hline

oa:6-13.poly & ${\frac {9}{80}}\,{x}^{6}+{\frac {43}{80}}\,{x}^{5}+{\frac {23}{16}}\,{x}^{4}+{\frac {143}{48}}\,{x}^{3}+{\frac {79}{20}}\,{x}^{2}+{\frac {179}{60}}
\,x+1$ \\ \hline

cut(4) &  ${\frac {2}{45}}\,{x}^{6}+{\frac {4}{15}}\,{x}^{5}+{\frac {7}{9}}\,{x}^{4}+
4/3\,{x}^{3}+{\frac {98}{45}}\,{x}^{2}+{\frac {12}{5}}\,x+1$ \\ \hline

cyclic01:5-8.poly & ${\frac {7}{60}}\,{x}^{5}+{\frac{5}{12}}\,{x}^{4}+5/4\,{x}^{3}+{\frac {31} {12}}\,{x}^{2}+{\frac{79}{30}}\,x+1$ \\ \hline

halfcube(5) & ${\frac {13}{15}}\,{x}^{5}+11/3\,{x}^{4}+16/3\,{x}^{3}+10/3\,{x}^{2}+9/5\,x
+1$ 
 \\ \hline
Cyclic(2,5)  & $10x^2+4x+1$ \\ \hline
Cyclic(3,5)  & $16x^3+10x^2+4x+1$ \\ \hline
Cyclic(4,5)  & $12x^4+16x^3+10x^2+4x+1$ \\ \hline
Cyclic(2,6)  & $20x^2+5x+1$ \\ \hline
Cyclic(3,6)  & $70x^3+20x^2+5x+1$ \\ \hline
Cyclic(4,6)  & $192x^4+70x^3+20x^2+5x+1$ \\ \hline
Cyclic(5,6)  & $288x^5+192x^4+70x^3+20x^2+5x+1$ \\ \hline
Cyclic(2,7)  & $35x^2+6x+1$ \\ \hline
Cyclic(3,7)  & $224x^3+35x^2+6x+1$ \\ \hline
Cyclic(4,7)  & $1512x^4+224x^3+35x^2+6x+1$ \\ \hline
Cyclic(2,8)  & $56x^2+7x+1$ \\ \hline
Cyclic(3,8)  & $588x^3+56x^2+7x+1$ \\ \hline
Cyclic(4,8)  & $8064x^4+588x^3+56x^2+7x+1$ \\ \hline
\end{tabular}
\caption{The Ehrhart polynomials for some well-known lattice polytopes. The choice
of coordinates for cyclic polytopes was $t=1,...,n$. The rest are listed Ehrhart
polynomials comes from $0/1$ polytopes selected from Ziegler's list. It includes
the Ehrhart polynomials of all $3$-dimensional $0/1$-polytopes.} 
\label{Ziegler}
\end{sidewaystable}
}

\subsection{$0/1$-polytopes}

\begin{figure}[htbp]
  \centering 
  \includegraphics[height=12cm]{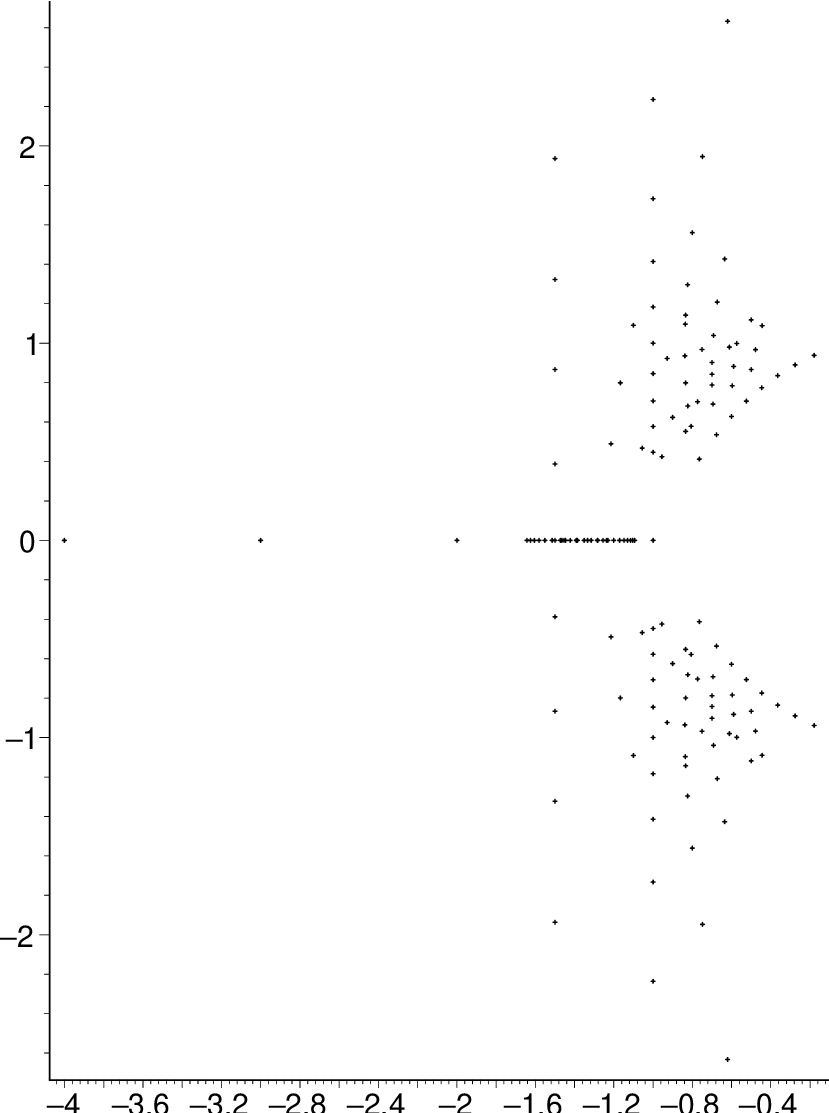}
  \includegraphics[height=12cm]{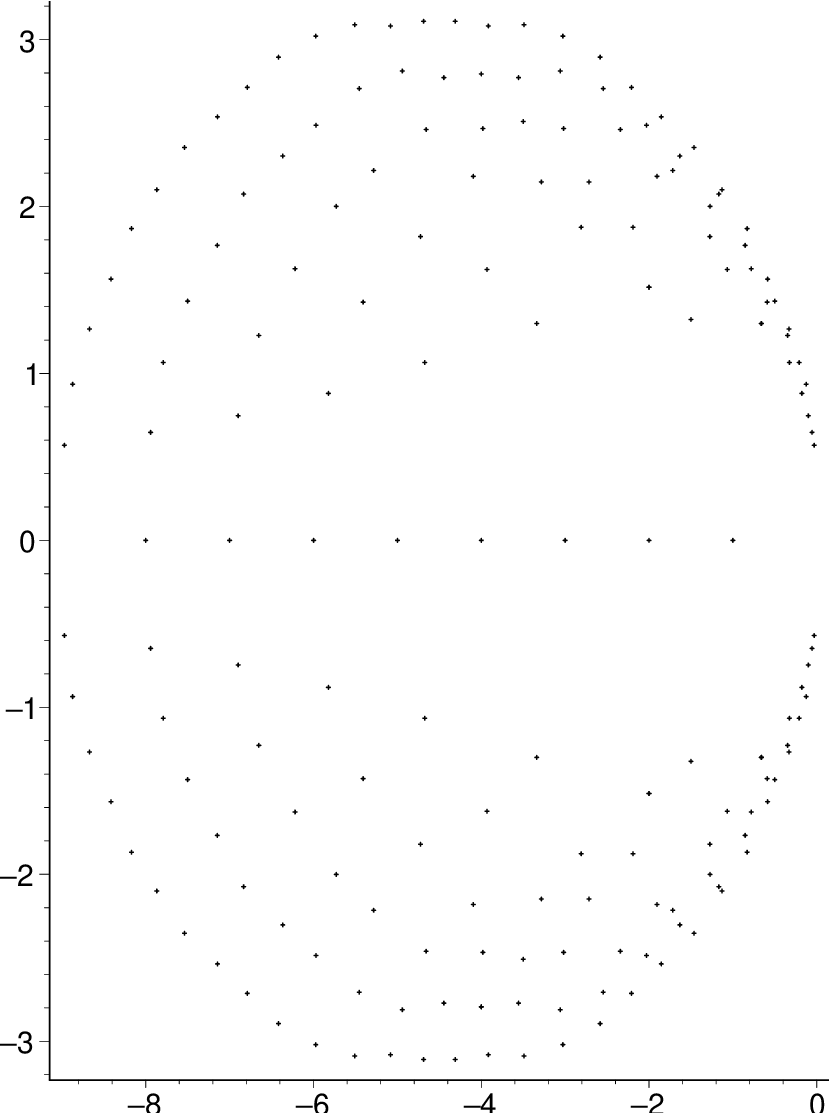}
 \caption{\emph{Top:} The zeros of the Ehrhart polynomials of all $3$ and $4$
dimensional $0/1$ polytopes. \emph{Bottom:} The zeros of Ehrhart polynomials
for the Birkhoff polytopes up to $n=9$}\label{fig:zerouno}
\end{figure}

We computed the Ehrhart polynomials for all $0/1$ polytopes of
dimension less or equal to 4 (up to symmetry there are 354 different
4-polytopes). In Figure \ref{fig:zerouno} we plotted their roots. In
our computations we relied on the on-line data sets of $0/1$ polytopes
available from {\tt Polymake}'s web page and those discussed in
Ziegler's lectures on $0/1$ polytopes~\cite{Kalai-Ziegler00}.  Several
phenomena are evident from the data we collected. For example, in
Table~\ref{Ziegler} we see two combinatorially different polytopes
that have the same Ehrhart polynomial. These are the so called ``nameless''
polytope of coordinates $(1, 0, 0, 0)$, $(1, 1, 0, 0)$, $(1, 0, 1,
0)$, $(1, 1, 1, 0)$, $(1, 0, 1, 1)$, $(1, 1, 0, 1)$ and the
octahedron. Another example of regular distribution appears also
in Figure \ref{fig:zerouno}. We show the roots of the
Ehrhart polynomials associated to the Birkhoff polytope of doubly
stochastic $n \times n$ matrices for $n=2, \dots, 9$.

\subsection{Cyclic polytopes}
Cyclic polytopes form a family whose combinatorial structure (i.e.
$f$-vector, face lattice, etc) is well understood. The canonical
choice of coordinates is given using the \emph {moment curve}
\begin{equation}
\label{eq:moment-curve} \nu_d: \left\{ \begin{array}{rcl} \mathbb{R} &
\to & \mathbb{R}^d,\\ t & \mapsto & (t^1,t^2,\dots,t^d) .
\end{array} \right.  \end{equation}
A cyclic polytope is obtained as the convex hull of $n$ points along
the moment curve. Thus we fix $t_1, t_2, \dots, t_n$ and define 
$C(n,d) := \conv \{ \nu_d(t_1), \nu_d(t_2),
\dots, \nu_d(t_n) \}$.  Cyclic polytopes are lattice polytopes exactly
when $t_i\in\Z$.  
There is a natural linear projection connecting these cyclic
polytopes.

\begin{lemma} 
  Consider the projection $\pi: \R^d\to\R^{d-1}$ that forgets the last
  coordinate. The inverse image under $\pi$ of a lattice point
  $y\in C(n,d-1)\cap\Z^{d-1}$ is a line that intersects the boundary of
  $C(n,d)$ in exactly two \emph{integral} points.
\end{lemma}

\begin{proof}
We need to prove that, given $t_1, t_2, \dots, t_d \in \Z$ and
$\lambda_1, \lambda_2, \dots, \lambda_d \in \R$,
\[
\forall \ 1 \leq j < d : \sum_{k=1}^d \lambda_k \, t_k^j \in \Z \qquad \Longrightarrow \qquad \sum_{k=1}^d \lambda_k \, t_k^d \in \Z \ .
\]
For $1 \leq j \leq d$, let $y_j= \sum_{k=1}^d \lambda_k \, t_k^j$; we
know that $y_1,y_2,\dots,y_{d-1} \in \Z$. We need to prove that
\[
\y = (y_1,\dots,y_d) = \sum_{k=1}^d \lambda_k \, \nu_d(t_k) \in \Z^d \ .
\]
This identity means that $\y$ lies on the hyperplane spanned by
$\nu_d(t_1),\dots,\nu_d(t_d)$, which can be expressed via a
determinant:
\[
\det \left( \begin{array}{cccc} 1 &       1 & \dots &       1 \\
                                  \y & \nu_d(t_1) &       & \nu_d(t_d) \end{array} \right) = 0 \ .
\]
Writing this determinant out through the first column and solving for $y_d$ gives
\begin{align*}
y_d \, = \, &- \frac 1 D \, \det \left( 
\begin{array}{ccc} t_1 & \cdots & t_d \\
                   \vdots & & \vdots \\
                   t_1^d & \cdots & t_d^d \end{array}
\right)
\, - \, \frac {y_1} D \, \det \left( 
\begin{array}{ccc} 1 & \cdots & 1 \\
                   t_1^2 & \cdots & t_d^2 \\
                   \vdots & & \vdots \\
                   t_1^d & \cdots & t_d^d \end{array}
\right)
\, - \, \cdots \\
&- \, \frac {y_{d-1}} D \, \det \left( 
\begin{array}{ccc} 1 & \cdots & 1 \\
                   t_1 & \cdots & t_d \\
                   \vdots & & \vdots \\
                   t_1^{d-2} & \cdots & t_d^{d-2} \\
                   t_1^d & \cdots & t_d^d \end{array}
\right) ,
\end{align*}
where 
\[
D \, = \, \det \left( \begin{array}{ccc} 1 & \cdots & 1 \\
                                t_1 & \cdots & t_d \\
                             \vdots &        & \vdots \\
                          t_1^{d-1} & \cdots & t_d^{d-1} \end{array} \right) = \prod_{ 1 \leq j < k \leq d } \left( t_j - t_k \right) .
\]
This expression yields an integer if we can prove that $D$ divides the
determinants appearing in the numerators. Equivalently, the
substitution $t_j = t_k$ in any of the numerators evaluates the
determinant to zero, which is apparent.
\end{proof}

Consequently, Conjecture~\ref{cyclicconj} is equivalent to saying that
the number of lattice points in a dilation of a cyclic polytope by a
positive integer $m$ is equal to its volume plus the number of lattice
points in its lower envelope. From the above lemma and Pick's theorem,
it follows that Conjecture~\ref{cyclicconj} is true for $d=2$.

\subsection*{Acknowledgements} We thank David Eisenbud, Francisco Santos, Bernd
Sturmfels, Tho\-mas Zaslavsky and G\"unter M. Ziegler for helpful
discussions and suggestions.  This research was supported in part by
the Mathematical Sciences Research Institute. Mike Develin was
also supported by the American Institute of Mathematics. Jes\'us De
Loera and Richard Stanley were partially supported by NSF grants
DMS-0309694 and DMS-9988459 respectively.

\bibliographystyle{plain} 
\bibliography{coeffzeros}
\setlength{\parskip}{0cm}

\end{document}

%% file: g-graphs.pstex_t
\begin{picture}(0,0)%
\includegraphics{g-graphs}%
\end{picture}%
\setlength{\unitlength}{4144sp}%
\begingroup\makeatletter\ifx\SetFigFont\undefined
\def\x#1#2#3#4#5#6#7\relax{\def\x{#1#2#3#4#5#6}}%
\expandafter\x\fmtname xxxxxx\relax \def\y{splain}%
\ifx\x\y   
\gdef\SetFigFont#1#2#3{%
  \ifnum #1<17\tiny\else \ifnum #1<20\small\else
  \ifnum #1<24\normalsize\else \ifnum #1<29\large\else
  \ifnum #1<34\Large\else \ifnum #1<41\LARGE\else
     \huge\fi\fi\fi\fi\fi\fi
  \csname #3\endcsname}%
\else
\gdef\SetFigFont#1#2#3{\begingroup
  \count@#1\relax \ifnum 25<\count@\count@25\fi
  \def\x{\endgroup\@setsize\SetFigFont{#2pt}}%
  \expandafter\x
    \csname \romannumeral\the\count@ pt\expandafter\endcsname
    \csname @\romannumeral\the\count@ pt\endcsname
  \csname #3\endcsname}%
\fi
\fi\endgroup
\begin{picture}(3428,2881)(870,-4484)
\put(2521,-3976){\makebox(0,0)[lb]{\smash{\SetFigFont{12}{14.4}{rm}{\color[rgb]{0,0,0}$B+1$}%
}}}
\put(2161,-4426){\makebox(0,0)[lb]{\smash{\SetFigFont{12}{14.4}{rm}{\color[rgb]{0,0,0}$\lambda s$ ($d$ even)}%
}}}
\put(3286,-4246){\makebox(0,0)[lb]{\smash{\SetFigFont{12}{14.4}{rm}{\color[rgb]{0,0,0}$\lambda s$ ($d$ odd)}%
}}}
\put(1261,-4021){\makebox(0,0)[lb]{\smash{\SetFigFont{12}{14.4}{rm}{\color[rgb]{0,0,0}$B$}%
}}}
\put(1216,-1906){\makebox(0,0)[lb]{\smash{\SetFigFont{12}{14.4}{rm}{\color[rgb]{0,0,0}$g$}%
}}}
\put(3916,-3391){\makebox(0,0)[lb]{\smash{\SetFigFont{12}{14.4}{rm}{\color[rgb]{0,0,0}$B+2$}%
}}}
\end{picture}